\documentclass[11pt]{amsart}

\usepackage{amsmath}
\usepackage{amssymb}
\usepackage{amscd}
\usepackage{color}
\usepackage{tikz}

\newcommand{\bccore}[1]{%
  \begin{tikzpicture}[baseline=(char.base)]
    \node[draw,circle,inner sep=0.1pt] (char) {\tiny #1};
  \end{tikzpicture}%
}

\begin{document}

\title{left dual (b,c)-core Inverses in Rings}

\author{Tugce Pekacar Calci}
\author{Serhat Emirhan Soycan}
\address{
Department of Mathematics, Ankara University, Ankara, Turkey}
\email{<tcalci@ankara.edu.tr>}
\address{
Ankara University Graduate School of National and Applied Sciences, Ankara, Turkey}
\email{<sesoycan@ankara.edu.tr>}

\newtheorem{thm}{Theorem}[section]
\newtheorem{lem}[thm]{Lemma}
\newtheorem{prop}[thm]{Proposition}
\newtheorem{cor}[thm]{Corollary}
\newtheorem{exs}[thm]{Examples}
\newtheorem{defn}[thm]{Definition}
\newtheorem{nota}{Notation}
\newtheorem{rem}[thm]{Remark}
\newtheorem{ex}[thm]{Example}
\newtheorem{que}[thm]{Question}


\begin{abstract} Let $a,b,c\in R$ where $R$ is a $*$-ring. We call $a$ \textit{left dual $(b,c)$-core invertible} if there exists $x\in Rc$ such that $bxab=b$ and $(xab)^*=xab$. Such an $x$ is called a left dual $(b,c)$-core inverse of $a$. In this paper, characteriztions of left dual $(b,c)$-core invertible element are introduced. We characterize left dual $(b,c)$-core inverses in terms of properties of the left annihilators and ideals. Moreover, we prove that $a$ is left dual $(b,c)$-core invertible if and only if $a$ is left $(b,c)$ invertible and $b$ is \{1,4\} invertible.  Also, properties of left dual $(b,c)$-core invertible elements are examined. We present the matrix representations of left dual $(b,c)$-core inverses by the Pierce decomposition. Furthermore, reletions between left dual $(b,c)$-core inverses and the other generalized inverses are given.
 \vspace{2mm}

\end{abstract}

\maketitle

\section{Introduction}
The Moore-Penrose inverse \cite{Penrose} and the Drazin inverse \cite{D1} are two important classes of generalized inverses. After those the inverse along an element \cite{Mary} and the $(b,c)$-inverse \cite{D2} which recover the Moore-Penrose and Drazin inverses were introduced. There are several researches about those inverses (see \cite{Bentez}, \cite{Boasso}, \cite{Ilic}, \cite{Marypat}). Shortly afterwards, one-sided inverses along an element were given \cite{Patricio}. Later, one-sided $(b,c)$-inverses were introduced by Drazin \cite{Drazin} which extends one-sided inverses along an element and $(b,c)$-inverses.

In 2023, Zhu introduced the $(b,c)$-core and dual $(b,c)$-core inverses in the context of *-semirgroups  which generalizes the core inverse \cite{Baksalary}, the core-EP inverse \cite{Mohana}, the Moore-Penrose inverse, the $w$-core inverse and right $w$-core inverse \cite{Mosic}.

Recently, authors introduced the right $(b,c)$-core inverses in a $*$-ring. 

Motivated by all of these, our aim is to introduce the left dual $(b,c)$-core inverses in a $*$-ring and give its relations with other classes of generalized inverses.  This offers a contribution to the theory of generalized inverses.  

In this paper, we introduce the notion of  left dual $(b,c)$-core inverse in a $*$-ring $R$. Let $a,b,c\in R$ where $R$ is a $*$-ring. We call $a$ \textit{left dual $(b,c)$-core invertible} if there exists $x\in Rc$ such that $bxab=b$ and $(xab)^*=xab$. Such an $x$ is called a left dual $(b,c)$-core inverse of $a$. Several fundamental results for characterizing left dual $(b,c)$-core inverses are established and it is shown that $a$ is left dual $(b,c)$-core invertible if and only if $a$ is left $(b,c)$ invertible and $b$ is \{1,4\} invertible. Moreover, $a_{l,\bccore{\#}\!(b,c)}=b^{(1,4)}a_l^{(b,c)}$. Then we characterize left dual $(b,c)$-core inverses in terms of the properties of the left annihilators and ideals. Further, we present the matrix representations of left dual $(b,c)$-core inverses by the Pierce decomposition.

After that, we give that several generalized inverses, such as left inverses, left dual core inverses, left dual pseudo core inverses, left dual $v$-core inverses, and Moore-Penrose inverses, are special situations of left dual $(b,c)$-core inverses. Precisely, for any nonnegative integers $m$ and $n$ satisfying $m + n = 1$, we establish the following equivalences in a ring:
 \begin{enumerate}
        \item $a$ is left invertible if and only if $a$ is left dual $(1,1)$-core invertible.
        \item $a$ is left dual pseudo core invertible if and only if $a^m$ is left dual $(a^k,a^n)$-core invertible for some positive integer $k$.
        \item $a$ is left dual core invertible if and only if $a^m$ is left dual $(a,a^n)$-core invertible.
        \item $a$ is left dual $v$-core invertible if and only if $v$ is left dual $(a,a)$-core invertible.
        \item $a$ is Moore-Penrose invertible if and only if $a$ is dual left $(a^*,a^*)$-core invertible if and only if $a^*$ is dual left $(a,a)$-core invertible.
    \end{enumerate}

In Section 3, we explore multiple characterizations of left dual core inverses and left dual $v$-core inverses. Additionally, we investigate the relationship between left dual pseudo core inverses and left dual $w$-core inverses.

Now, let's remember a few concepts about generalized inverses. Let $R$ be an associative ring with unity $1$. An element $a \in R$ is called \textit{(von Neumann) regular} if there exists some $x \in R$ such that $axa = a$. Such an $x$ is called an \textit{inner inverse} or a $\{1\}$-inverse of $a$ and is denoted by $a^-$. The symbol $a\{1\}$ means the set of all inner inverses of $a$. The set of all regular elements in $R$ is denoted by $R^-$.

In \cite{Patricio}, Zhu et al. extended inverses along an element to one-sided cases. Let $a, d \in R$. An element $a$ is called \textit{left invertible along} $d$ if there exists some $x \in R$ such that $xad = d$ and $x \in Rd$. Such an element $x$ is called a \textit{left inverse of $a$ along $d$} and is denoted by $a_l^d$. Dually, an element $a$ is called right invertible along $d$ if there exists some $y\in R$ such that $day=d$ and $y\in dR$. Such element $y$ is called a right inverse of $a$ along $d$ and is indicated by $a_r^d$.

In 2016, Drazin defined one-sided \((b, c)\)-inverses. Let \(a, b, c \in R\). We call \(a\) left \((b, c)\)-invertible if \(b \in Rcab\), or equivalently if there exists \(x \in Rc\) such that \(xab = b\). In this case, any such \(x\) will be called a left \((b, c)\)-inverse of \(a\) and denoted by \(a_{l}^{(b,c)}\). 

Dually, \(a\) is right \((b, c)\)-invertible if \(c \in cabR\), or equivalently if there exists \(y \in bR\) such that \(cay = c\). In this case, any such \(y\) will be called a right \((b, c)\)-inverse of \(a\) and denoted by \(a_{r}^{(bc)}\). In particular, \(a\) is called \((b, c)\)-invertible if it is both left and right \((b, c)\)-invertible. 

We denote by \(R_{l}^{(b,c)}\), \(R_{r}^{(b,c)}\) and \(R^{(b,c)}\), the set of all left \((b, c)\)-invertible, the set of all right \((b, c)\)-invertible and the set of all \((b, c)\)-invertible elements in \(R\), respectively. It should be noted that \(a\) is left \((d, d)\)-invertible if and only if it is left invertible along \(d\). Moreover, the left \((d, d)\)-inverse of \(a\) is exactly the left inverse of \(a\) along \(d\).

A map \(\ast : R \to R\) is an \textit{involution of \(R\)} if it satisfies \((x^\ast)^\ast = x\), \((xy)^\ast = y^\ast x^\ast\) and \((x+y)^\ast = x^\ast + y^\ast\) for all \(x, y \in R\). Throughout this paper, any ring \(R\) is assumed to be a unital ring, that is, a ring \(R\) with unity \(1\) and an involution \(\ast\).

An element \(a \in R\) is said to be \textit{Moore-Penrose invertible} if there exists some \(x \in R\) such that $axa = a,~xax = x,~ (ax)^\ast = ax~ \text{and}~ (xa)^\ast = xa.$
Such an \(x\) is called a Moore-Penrose inverse of \(a\). It is unique if it exists and is denoted by \(a^\dagger\). Generally, any solution \(x\) satisfying the equations \(axa = a\) and \((ax)^\ast = ax\) (respectively, \((xa)^\ast = xa\)) is called a \{1,3\}-inverse (respectively, a \{1,4\}-inverse) of \(a\). 

The symbols \(a^{(1,3)}\) and \(a^{(1,4)}\) denote a \{1,3\}-inverse and a \{1,4\}-inverse of \(a\), respectively. We denote by \(a\{1,3\}\) and \(a\{1,4\}\), the set of all \{1,3\}-inverses and \{1,4\}-inverses of \(a\), respectively. In general, the set of all \{1,3\}-invertible, \{1,4\}-invertible and Moore-Penrose invertible elements in \(R\) will be denoted by \(R^{\{1,3\}}\), \(R^{\{1,4\}}\) and \(R^{\dagger}\), respectively. It is known that \(a\) is Moore-Penrose invertible if and only if it is both \{1,3\}-invertible and \{1,4\}-invertible. An element \(p \in R\) is called a \textit{projection} if \(p^2 = p = p^\ast\).

An element \(a \in R\) is called \textit{left dual pseudo core invertible} if there exist \(x \in R\) and a positive integer \(k\) such that $a^kxa=a^k$, $(xa)^*=xa$ and $x^2a=x$. Such an $x$ is called a left dual pseudo core inverse of $a$ and denoted by $a_{l,\bccore{D}}$. In this case, the smallest positive integer $k$ is called left dual pseudo core index of $a$ and denoted by $I(a)$. In particular, $a$ is called left dual core invertible if $I(a)=1$. Generally, $R_{l,\bccore{D}}$ and $R_{l,\bccore{\#}}$ denote the set of all left dual pseudo core invertible and left dual core invertible elements in $R$, respectively.

For any $a,w\in R$, $a$ is called \textit{left dual $v$-core invertible} if there exists some $x\in R$ such that $axva=a$, $xva=(xva)^*$ and $x^2va=x$. Such an $x$ is called a left dual $v$-core inverse of $a$ and denoted by $a_{l,v\bccore{\#}}$. The symbol $R_{l,v\bccore{\#}}$ denotes the set of all left dual $v$-core invertible elements in $R$.

The \((b, c)\)-core inverse was defined in an \(a\)-monoid \(S\) . For convenience, we next state this notion in \(R\). Let \(a, b, c \in R\). The element \(a\) is called \textit{\((b, c)\)-core invertible} if there exists some \(x \in R\) such that:
\[
caxc = c, \quad xR = bR \quad \text{and} \quad Rx = Rc^*.
\]

Dually, an element $a\in R$ is called \textit{dual $(b,c)$-core invertibl}e if there exists some $y\in R$ such that: $$byab = b, \quad yR = b^*R \quad \text{and} \quad Ry = Rc.$$
The $(b, c)$-core inverse and dual $(b, c)$-core inverse of $a$ are unique if they exist and are denoted by $a_{(b,c)}^{\bccore{\#}\!}$ and $a_{\bccore{\#}\!(b,c)}$, respectively. As usual, we denote by $R_{(b,c)}^{\bccore{\#}\!}$ and $R_{\bccore{\#}\!(b,c)}$, the set of all \((b, c)\)-core invertible and dual \((b, c)\)-core invertible elements in \(R\), respectively.

\section{Left Dual (b,c)-core inverses}
In this section, we define left dual $(b,c)$-core invertible elements and we explore multiple characterizations of them.
\begin{defn}
    Let $a,b,c\in R$. We call $a$ left dual $(b,c)$-core invertible if there exists some $x\in Rc$ such that $bxab=b$ and $(xab)^*=xab$. Such an $x$ is called a left dual $(b,c)$-core inverse of $a$.  The set of all left dual $(b,c)$-core invertible elements in $R$ is denoted by $R_{l,\bccore{\#}\!(b,c)}$.
\end{defn}

We will denote the left dual $(b,c)$-core inverse of $a$ with $a_{l,\bccore{\#}\!(b,c)}$. An element $a$ in $R$ could have different left dual $(b,c)$-core inverses. For example, let $R$ be a $*$-ring, $a\in R$, $b=0$ and $c=1$. Then every $x\in R$ is the left dual $(0,1)$-core inverse of $a$.

It is obvious that every dual $(b,c)$-core invertible element is left dual $(b,c)$-core invertible. The converse statement is not always true. For example, let $R$ be a $*$-ring, $a\in R$, $b=0 \neq c$. Clearly, $a$ is left dual $(b,c)$-core invertible. But $a$ is not dual $(b,c)$-core invertible since  $cabx=0$ for every $x\in R$.
\begin{thm}
    Let $R$ be a ring and $a,b,c\in R$. The following conditions are equivalent.
    \begin{enumerate}
        \item $a\in R_{l,\bccore{\#}\!(b,c)}$.
        \item There exists some $x\in Rc$ such that $bxab=b$, $(xab)^*=xab$ and $xabx=x$.
        \item There exists some $x\in R$ such that $bxab=b$, $xR=b^*R$ and $Rx\subseteq Rc$.
        \item There exists some $x\in R$ such that $bxab=b$, $l(x)=l(b^*)$ and $Rx\subseteq Rc$.
        \item There exists some $x\in R$ such that $bxab=b$, $Rx\subseteq Rc$ and $xR\subseteq b^*R$.
        \item There exists some $x\in R$ such that $bxab=b$, $Rx\subseteq Rc$ and  $l(b^*)\subseteq l(x)$.
        \item There exists a projection $q\in R$ and an idempotent $p \in R$ such that $Rb\subseteq Rq\subseteq Rab$, $Rp\subseteq Rc$ and $abR\subseteq pR$.
    \end{enumerate}
    In this case, $a_{l,\bccore{\#}\!(b,c)}=q(ab)^-p$ for any $(ab)^-\in (ab)\{1\}$.
\end{thm}

\begin{proof}
    (1) $\Rightarrow$ (2) Assume $a\in  R_{l,\bccore{\#}\!(b,c)}$. Then there exists some $y\in Rc$ such that $byab=b$ and $(yab)^*=yab$. Let $x=yaby$. We get $xab=(yaby)ab=ya(byab)=yab$. And so $(xab)^*=xab$. $bxab=b(yaby)ab=(byab)yab=byab=b$. $xabx=(yaby)ab(yaby)=ya(byab)yaby=ya(byab)y=yaby=x$.

    (2) $\Rightarrow$ (3) Since $bxab=b$ and $(xab)^*=xab$, we have $b^*=xabb^*\in xR$. Also, $xabx=x$ implies $x=(xab)^*x=b^*a^*x^*x$. Therefore, $x\in b^*R$.
    
    (3) $\Rightarrow$ (4) Obvious.

    (4) $\Rightarrow$ (5) Since $b=bxab$, $b^*=b^*a^*x^*b^*$. Therefore $b^*-b^*a^*x^*b^*=(1-b^*a^*x^*)b^*=0$. Since $l(x)=l(b^*)$, $(1-b^*a^*x^*)x=0$. And so $x=b^*a^*x^*x\in b^*R$, as desired.

    (5) $\Rightarrow$ (6) Obvious.

    (6) $\Rightarrow$ (7) Since $b^*=(xab)^*b^*$ and $l(b^*)\subseteq l(x)$, we get $x=(xab)^*x$. Then we have $xab=(xab)^*(xab)$ which implies $(xab)^*=(xab)$. Set $q=xab$ and $p=abx$. Hence $q^2=(xab)(xab)=xa(bxab)=xab$ and $q^*=(xab)^*=xab=q$. So $q$ is a projection. Also $p^2=(abx)(abx)=a(bxab)x=abx=p$, so $p$ is an idempotent. Therefore, $Rb=Rbq\subseteq Rq=R(xab)\subseteq Rab$, $Rp=R(abx)\subseteq Rx\subseteq Rc$ and $abR=a(bxab)R=(abx)abR=pabR\subseteq pR$.

    (7) $\Rightarrow$ (1) Since $abR\subseteq pR$, $pab=ab$. From $Rb\subseteq Rq\subseteq Rab$, it follows that $b=bq$ and $q=zab$ for some $z\in R$. Therefore, $ab=a(bq)=ab(zab)=(ab)z(ab)$ and so $ab\in R^-$. Let $x=q(ab)^-p$ for any $(ab)^-\in (ab)\{1\}$. Then, $x=q(ab)^-p\in Rp\subseteq Rc$. Also, $xab=q(ab)^-pab=q(ab)^-ab=zab(ab)^-(ab)=zab=q=(xab)^*$. Moreover, we have $bxab=bq=b$, as asserted.
\end{proof}

\begin{lem}\cite[Lemma 2.2]{Zhang} \label{zhang} 
    Let $R$ be a ring and $a \in R$. Then $a^{(1,4)}$ exists if and only if $a\in aa^*R$. If $aa^*y=a$ for some $y\in R$, then $y^*=a^{(1,4)}$.
\end{lem}

Assume that $a\in R_{l,\bccore{\#}\!(b,c)}$ and a left dual $(b,c)$-core inverse of $a$ is $x$. Then since $bxab=b$, we get $(xab)^n=xab$ for any positive integer $n$. It is concluded that $a\in R_{l,\bccore{\#}\!(b,c)}$ implies $x\in Rc$, $b(xab)^n=b$ and $((xab)^n)^*=(xab)^n$ for any positive integer $n$. It is natural to ask if the converse implication is true. The following theorem provides a solution to that question.

\begin{thm}\label{denklik}
    Let $R$ be a ring and $a,b,c\in R$. The following statements are equivalent.
    \begin{enumerate}
        \item $a\in R_{l,\bccore{\#}\!(b,c)}$.
        \item $b\in b(cab)^*R$.
        \item $b\in Rcab\cap bb^*R$.
        \item There exists some $x\in Rc$ such that $b(xab)^n=b$ and $((xab)^n)^*=(xab)^n$ for any positive integer $n$.
        \item There exists some $x\in Rc$ such that $b(xab)^n=b$ and $((xab)^n)^*=(xab)^n$ for some positive integer $n$. 
    \end{enumerate}
\end{thm}

\begin{proof}
    (1) $\Rightarrow$ (2) Suppose $a\in  R_{l,\bccore{\#}\!(b,c)}$. Then there exists some $x\in Rc$ such that $bxab=b$, $(xab)^*=xab$. Let $x=rc$. Then $b=bxab=b(xab)^*=b(rcab)^*=b(cab)^*r^*\in b(cab)^*R$, as desired.

    (2) $\Rightarrow$ (3) It is clear by \cite[Lemma 3.1 II]{Zhu}.

    (3) $\Rightarrow$ (4) Since $b\in Rcab\cap bb^*R$, we have $b=rcab=bb^*s$ for some $r,s\in R$. In this case, $s^*=b^{(1,4)}$ by Lemma \ref{zhang}. Let $x=s^*rc$. Then $x\in Rc$, $xab=(s^*rc)ab=s^*(rcab)=s^*b=(s^*b)^*=(xab)^*$ and $bxab=bs^*b=b$. Multiplying the latter equation by $xa$ from left, we get $(xab)^2=xab$, and so $(xab)^n=xab$ for any positive integer $n$. As a consequence $b(xab)^n=bxab$ and $((xab)^n)^*=(xab)^*=xab=(xab)^n$.

    (4) $\Rightarrow$ (5) Obvious.
    
    (5) $\Rightarrow$ (1) Suppose $x\in Rc$ such that $b(xab)^n=b$ and $((xab)^n)^*=(xab)^n$ for some positive integer $n$. Set $y=(xab)^{(n-1)}x$. We claim that $y=a_{l\bccore{\#}\!(b,c)}$. Indeed, $y=(xab)^{(n-1)}x \in Rc$. Also, $byab=b((xab)^{(n-1)}x)ab=b(xab)^n=b$. Moreover, $yab=((xab)^{(n-1)}x)ab=(xab)^n=((xab)^n)^*=(yab)^*$.
\end{proof}

Drazin \cite{Drazin} introduced the concept of left and right $(b,c)$-invertible elements. Now we shall give this definition for completeness.

\begin{defn}
    Let $S$ be any semigroup and $a,b,c\in S$. Then we shall say that $a$ is left $(b,c)$-invertible if $b\in Scab$, or equivalently if there exists $x\in Sc$ such that $xab=b$, in which case any such $x$ will be called a left $(b,c)$-inverse of $a$.

    Dually $a$ is called right $(b,c)$-invertible if $c\in cabS$, or equivalently if there exists $z\in bS$ such that $caz=c$, and any such $z$ will be called a right $(b,c)$-inverse of $a$.
\end{defn}
Let $a,b,c \in R$. Zhu showed that $a$ is dual $(b,c)$-core invertible if and only if $a$ is $(b,c)$-invertible and $b$ ($ab$ or $cab$) is \{1,4\}-invertible. An analogous result on dual right $(b,c)$-core inverses can be obtained.

\begin{thm}\label{hehehe}
    Let $R$ be a ring and $a,b,c\in R$. The following statements are equivalent.
    \begin{enumerate}
        \item $a\in R_{l,\bccore{\#}\!(b,c)}$.
        \item $a\in R_l^{(b,c)}$ and $b \in R^{(1,4)}$.
        \item $a\in R_l^{(b,c)}$ and $ab \in R^{(1,4)}$.
        \item $a\in R_l^{(b,c)}$ and $cab \in R^{(1,4)}$.
        \end{enumerate}
        \noindent In this case $a_{l,\bccore{\#}\!(b,c)}=b^{(1,4)}a_l^{(b,c)}=(ab)^{(1,4)}aa_l^{(b,c)}=(cab)^{(1,4)}c$.
\end{thm}

\begin{proof}
   (1) $\Rightarrow$ (2) It is a direct consequence of Lemma \ref{zhang} and Theorem \ref{denklik}.

   (2) $\Rightarrow$ (3) Let $a\in R_l^{(b,c)}$ and $b\in R^{(1,4)}$. Then there exists some $s\in R$ such that $scab=b$. Let $x=b^{(1,4)}$. Now, $abxscab=abxb=ab$. So $ab$ is regular and $xscab=xb=(xb)^*=(xscab)^*$. Hence, $ab \in R^{(1,4)}$.

   (3) $\Rightarrow$ (4) Let $a\in R_l^{(b,c)}$ and $ab\in R^{(1,4)}$. Then there exists some $s\in R$ such that $scab=b$. Suppose that $x=ab^{(1,4)}$. Now $cabxascab=cabxab=cab$. So $cab$ is regular and $xascab=xab=(xab)^*=xascab$. Hence $cab \in R^{(1,4)}$.

   (4) $\Rightarrow$ (1)  Let $a\in R_l^{(b,c)}$ and $cab\in R^{(1,4)}$. Then there exists $s\in R$ such that $b=scab$. Let $x=cab^{(1,4)}$. Set $y=xc$. Then $y\in Rc$, $yab=xcab=(xcab)^*=(ycab)^*$ and $byab=bxcab=scabxcab=scab=b$. Therefore, $a\in  R_{l,\bccore{\#}\!(b,c)}$ . The other equalities can be seen in a similar way.
\end{proof}

Suppose that $a\in R_{l,\bccore{\#}\!(b,c)}$. Theorem \ref{hehehe} guarantees that $cab\in R^{(1,3)}$. Therefore, $cab\in R^-$. Also by Theorem \ref{denklik}, there exists some $s\in R$ such that $b=scab$. It follows that $b=scab(cab)^-cab=b(cab)^-cab$ for any $(cab)^-\in (cab)\{1\}$. This implies $a_l^{(b,c)}=b(cab)^-c$. Now another representation of $a_{l,\bccore{\#}\!(b,c)}$ can be presented.

\begin{prop}
    Let $R$ be a ring and $a,b,c\in R$ with $a\in R_{l,\bccore{\#}\!(b,c)}$. Then $a_{l,\bccore{\#}\!(b,c)}=b^{(1,4)}b(cab)^-c$ for any $(cab)^-\in (cab)\{1\}$ and $b^{(1,4)} \in b\{1,4\}$.
\end{prop}

\begin{thm}
    Let $R$ be a ring and $a,b,c\in R$. Then the following statements are equivalent.
    \begin{enumerate}
        \item $a\in R_{l,\bccore{\#}\!(b,c)}$.
        \item $R=(cab)^*R\oplus r(b)$.
        \item $R=(cab)^*R+r(b)$.
    \end{enumerate}
\end{thm}

\begin{proof}
(1) $\Rightarrow$ (2) Let $a\in R_{l,\bccore{\#}\!(b,c)}$. By Theorem \ref{denklik}, we have $b\in b(cab)^*R$. Hence $b=bb^*a^*c^*r$ for some $r\in R$. So $1-b^*a^*c^*r=1-(cab)^*r \in r(b)$. For any $s\in R$, we have $s=[(1-(cab)^*r)+(cab)^*r]s= (1-(cab)^*)s+(cab)^*rs\in (1-(cab)^*)R+(cab)^*R$. Since $b=bb^*a^*c^*r$ and Lemma \ref{zhang}, $r^*ca \in b\{1,4\}$. Then for any $z\in (cab)^*R \cap r(b)$, $bz=0$ and there exists some $s\in R$ such that $z=(cab)^*s=(cab(r^*ca)b)^*s=(r^*cab)^*(cab)^*s=r^*cabz=0$. Hence, $R=(cab)^*R \oplus r(b)$.

    (2) $\Rightarrow$ (3) Obvious.

    (3) $\Rightarrow$ (1) Let $R=(cab)^*R+r(b)$. Then $b\in b(cab)^*R$. By Theorem \ref{denklik}, $a\in R_{l,\bccore{\#}\!(b,c)}$ as desired.
\end{proof}

For any $p\in R$, any $a\in R$ can be written as 
$$a=pap+pa(1-p)+(1-p)ap+(1-p)a(1-p) $$ or the matrix form 

$$a=\left[%
\begin{array}{cc}
  a_1 & a_2 \\
  a_3 & a_4 \\
\end{array}%
\right]_p ,$$

\noindent where $a_1=pap$, $a_2=pa(1-p)$, $a_3=(1-p)ap$ and $a_4=(1-p)a(1-p)$. This decomposition is well known as the Pierce decomposition.

If $p^2=p=p^*$, then

$$a^*=\left[%
\begin{array}{cc}
  {a_1}^* & {a_2}^* \\
  {a_3}^* & {a_4}^* \\
\end{array}%
\right]_p.$$

We next give the matrix represantations of left dual $(b,c)$-core inverses.

\begin{thm}
    Let $R$ be a ring, $a,b,c\in R$ and $\textbf{R}(c)$ denote the column space of $c$. Then the following statements are equivalent.
    \begin{enumerate}
        \item $a\in R_{l,\bccore{\#}\!(b,c)}$ and $x$ is a left dual $(b,c)$-core inverse of $a$.
        \item There exists a projection $q\in R$ such that 
         $$a=\left[%
\begin{array}{cc}
  {a_1} & {a_2} \\
  {a_3} & {a_4} \\
\end{array}%
\right]_q,~b=\left[%
\begin{array}{cc}
  {b_1} & 0 \\
  {b_3} & 0 \\
\end{array}%
\right]_q,~ c=\left[%
\begin{array}{cc}
  {c_1} & {c_2} \\
  {c_3} & {c_4} \\
\end{array}%
\right]_p\text{and } x=\left[%
\begin{array}{cc}
  {x_1} & {x_2} \\
  {x_3} & {x_4} \\
\end{array}%
\right]_q$$ 
\noindent where $(x_1a_1+x_2a_3)b_1+(x_1a_2+x_2a_4)b_3=q$, $(x_3a_1+x_4a_3)b_1+(x_3a_2+x_4a_4)b_3=0$ and $\textbf{R}(x)\subseteq \textbf{R}(c)$.
\item  There exists a projection $p\in R$ such that 
 $$ a=\left[%
\begin{array}{cc}
  {a_1} & {a_2} \\
  {a_3} & {a_4} \\
\end{array}%
\right]_p,~b=\left[%
\begin{array}{cc}
  0 & {b_2} \\
  0 & {b_4} \\
\end{array}%
\right]_p,~ c=\left[%
\begin{array}{cc}
  {c_1} & {c_2} \\
  {c_3} & {c_4} \\
\end{array}%
\right]_p\text{ and } x=\left[%
\begin{array}{cc}
  {x_1} & {x_2} \\
  {x_3} & {x_4} \\
\end{array}%
\right]_p$$
where $(x_1a_1+x_2a_3)b_2+(x_1a_2+x_2a_4)b_4=1-p$, $(x_3a_1+x_4a_3)b_2+(x_3a_2+x_4a_4)b_4=0$ and $\textbf{R}(x)\subseteq \textbf{R}(c)$.
    \end{enumerate}
\end{thm}

\begin{proof}
    (1) $\Rightarrow$ (2) Suppose $a\in R_{l,\bccore{\#}\!(b,c)}$ with a left dual $(b,c)$-core inverse $x$. Then $x\in Rc$, $(xab)^*=xab$ and $bxab=b$. Let $q=xab$. Then $q^2=q=q^*$. So, $a,b,c$ can be represented as in $(2)$. Since $x\in Rc$, $\textbf{R}(x)\subseteq \textbf{R}(c)$. By the Pierce decomposition, we have 
    $$(x_1a_1+x_2a_3)b_1+(x_1a_2+x_2a_4)b_3$$
    $$=(qxqqaq+qx(1-q)(1-q)aq)qbq+(qxqqa(1-q)+qx(1-q)(1-q)a(1-q))(1-q)bq$$
    $$=qxqaqb+qxaqb-qxqaqb+qxqa(1-q)b+qxa(1-q)b-qxqa(1-q)b $$
    $$\hspace{-5.05cm}=qxaqb+qxqab-qxqaqb+qxab-qxaqb$$
    $$\hspace{-10.8cm}=qxab$$
    $$\hspace{-11.4cm}=q.$$

    The equality $(x_3a_1+x_4a_3)b_1+(x_3a_2+x_4a_4)b_3=0$ can be shown in a similar way.

    (2) $\Rightarrow$ (1)
   Since $$xab=\left[%
\begin{array}{cc}
  {(x_1a_1+x_2a_3)b_1+(x_1a_2+x_2a_4)b_3} & 0 \\
  {(x_3a_1+x_4a_3)b_1+(x_3a_2+x_4a_4)b_3} & 0 \\
\end{array}%
\right]_q=\left[%
\begin{array}{cc}
  q & 0 \\
  0 & 0 \\
\end{array}%
\right]_q=q, $$ here it is obvious that $bxab=b$ and $(xab)^*=xab$. Besides, since $\textbf{R}(x)\subseteq \mathbf{R}(c)$, $x\in Rc$. Hence, $a\in R_{l,\bccore{\#}\!(b,c)}$ and $x$ is a left dual $(b,c)$-core inverse of $a$.

(1) $\Leftrightarrow$ (3) It can be shown in a similar way of the proof of (1) $\Leftrightarrow$ (2) by taking $p=1-q$.
\end{proof}

\section{Applications}

In this section, we explore multiple characterizations of left dual core inverses and left dual $v$-core inverses. Additionally, we investigate the relationship between left dual pseudo core inverses, left dual $w$-core inverses and left dual $(b,c)$-core inverses.

\begin{thm}\label{relation}
    Let $R$ be a ring and $a,w \in R$ and $m,n$ be nonnegative integers such that $m+n\geq 1$. Then 
    \begin{enumerate}
        \item $a$ is left invertible if and only if $a$ is left dual $(1,1)$-core invertible,
        \item $a$ is left dual pseudo core invertible if and only if $a^m$ is left dual $(a^k,a^n)$-core invertible, for some positive integer $k$,
        \item $a$ is left dual core invertible if and only if $a^m$ is left dual $(a,a^n)$-core invertible,
        \item $a$ is left dual $v$-core invertible if and only if $v$ is left dual $(a,a)$-core invertible.
    \end{enumerate}
\end{thm}

\begin{proof}
    (1) Let $a$ be left dual $(1,1)$-core invertible. Then there exists some $x\in R$ such that $(xa)^*=xa$ and $xa=1$. Converse is clear.

    (2) Let $a$ be left dual pseudo core invertible with $I(a)=k$. Then there exists some $x\in R$ such that $a^kxa=a^k$, $(xa)^*=xa$ and $x^2a=x$, whence $xa=x^2aa=x^2a^2=\cdots=x^na^n$ for any positive integer $n$. Let $y=x^{k+m}$. Then
    $$y=x^{k+m}=x^{k+m-1}x^2a=x^{k+m-1}x^{n+1}a^n=x^{k+m+n}a^n\in Ra^n,$$
    $$ya^ma^k=ya^{m+k}=x^{m+k}a^{m+k}=xa=(xa)^*=(ya^ma^k)^*,$$ 
          $$a^kya^ma^k=a^kx^{k+m}a^{m+k}=a^kxa=a^k.$$
    Hence, $a^m\in R_{l,\bccore{\#}\!(a^k,a^n)}$ and $(a^m)_{l,\bccore{\#}\!(a^k,a^n)}=(a_{l,\bccore{D}})^{k+m}$. 

    Conversely, let $a^m\in R_{l\bccore{\#}\!(a^k,a^n)}$. Then there exists some $x\in Ra^n$ such that $(xa^ma^k)^*=(xa^{k+m})^*=xa^{k+m}$ and $a^kxa^ma^k=a^k$. Let $z=xa^{k+m-1}$. Then
     $a^kza=a^kxa^{k+m-1}a=a^kxa^{k+m}=a^kxa^ma^k=a^k$. Also, $za=xa^{k+m-1}a=xa^{k+m}=(xa^{k+m})^*=(za)^*$. Moreover, we have  $z^2a=xa^{k+m-1}xa^{k+m-1}a=xa^{k+m-1}xa^{k+m}=xa^{m-1}a^kxa^ma^k=xa^{m-1}a^k=xa^{k+m-1}=z$. Therefore $a\in R_{l\bccore{D}}$ and $a_{l,\bccore{D}}=(a^m)_{l,\bccore{\#}\!(a^k,a^n)}a^{k+m-1}$.
    
    (3) It is obvious from (2).

(4) Let $a\in  R_{l,v,\bccore{\#}}$. Then there exists some $y\in R$ such that $ayva=a$, $yva=(yva)^*$ and $y^2va=y$. Let $x=y$. Then,  $x=x^2va\in Ra$ and $(xva)^*=xva$. Also, $axva=a$. $v\in R_{l,\bccore{\#}\!(a,a)}$ and $v_{l,\bccore{\#}\!(a^k,a^n)}=a_{l,v,\bccore{\#}}$. Now let $v\in R_{l,\bccore{\#}\!(a,a)}$. Then there exists some $x\in Ra$ such that $(xva)^*=(xva)$ and $axva=a$. Let $y=x$, then there exists some $r\in R$ such that $y=ra$. Hence, $ayva=a$ and $(yva)^*=yva$. Moreover, we have $y^2va=yyva=rayva=ra=y$, as asserted.
 \end{proof}

It is proved in Theorem \ref{relation} that $a$ is left dual $v$-core invertible if and only if $v$ is left dual $(a,a)$-core invertible. As a result of Theorem \ref{relation}, we get the matrix representation of the left dual $v$-core inverses as follows.

\begin{cor}
    Let $R$ be a ring and $a,v\in R$. Then the following statements are equivalent:

    \begin{enumerate}
        \item $a\in R_{l,v,\bccore{\#}}$ and $x$ is a left dual $v$-core inverse of $a$.
        \item There exists a projection $q\in R$ such that 
        
        $v=\left[%
\begin{array}{cc}
  v_1 & v_2 \\
  v_3 & v_4 \\
\end{array}%
\right]_q$, $a=\left[%
\begin{array}{cc}
  a_1 & 0 \\
  a_3 & 0 \\
\end{array}%
\right]_q$ $x=\left[%
\begin{array}{cc}
  x_1 & 0 \\
  x_3 & 0 \\
\end{array}%
\right]_q$

 where $x_1v_1a_1+x_1v_2a_3=q$ and $x_3v_1a_1+x_3v_2a_3=0$.
 \item  There exists a projection $q\in R$ such that 
        
        $v=\left[%
\begin{array}{cc}
  v_1 & v_2 \\
  v_3 & v_4 \\
\end{array}%
\right]_p$, $a=\left[%
\begin{array}{cc}
  0 & a_2 \\
  0 & a_4 \\
\end{array}%
\right]_p$ $x=\left[%
\begin{array}{cc}
  0 & x_2 \\
  0 & x_4 \\
\end{array}%
\right]_p$

 where $x_2v_3a_2+x_2v_4a_4=1-p$ and $x_4v_3a_2+x_4v_4a_4=0$.
    \end{enumerate}
   
\end{cor}

\begin{thm}\label{dec}
    Let $R$ be a ring and $a,w\in R$ with $a\in R_{l,v,\bccore{\#}}$. Then $va=a_1+a_2$ where
    \begin{enumerate}
        \item $a_1\in R_{l,\bccore{\#}}$,
        \item ${a_2}^2=0$,
        \item ${a_2}^*a_1=0=a_1a_2.$
        In addition, $a_{l,v,\bccore{\#}}(va)^2 \in R_{l,\bccore{\#}}$ with a left dual core inverse $a_{l,v,\bccore{\#}}$.
    \end{enumerate}
\end{thm}

\begin{proof}
    Suppose $a\in R_{l,v,\bccore{\#}}$ with a left dual $v$-core inverse $x$. Then $axva=a$, $(xva)^*=xva$ and $x^2va=x$. Let $a_1=x(va)^2$ and $a_2=(1-xva)va$. Then $va=a_1+a_2$. Now,

    (1) $xa_1=xx(va)^2=x^2(va)(va)=xva=(xva)^*=(xa_1)^*$. Also, we have $a_1xa_1=a_1(xva)=x(va)^2xva=x(va)(va)xva=xvava=x(va)^2=a_1$. Moreover, $x^2a_1=xxa_1=x(xva)=x^2va=x$. Hence, $a_1\in R_{l,\bccore{\#}}$ with a left dual core inverse $a_{l,\bccore{\#}}$. 
    
    (2) ${a_2}^2=(va-xvava)^2=vava-vax(va)^2-x(va)^3+x(va)^2x(va)^2=vava-vava-x(va)^3+x(va)^3=0$. 
    
    (3) ${a_2}^*a_1=(va-x(va)^2)^*x(va)^2=((va)^*-(va)^*xva)x(va)^2=(va)^*(1-xva)x(va)^2=(va)^*(x(va)^2-xv(axva)va)=(va)^*(x(va)^2-x(va)^2)=0$. $a_1a_2=0$ can be shown in a similar way.
 \end{proof}
Let $a\in R$ and $m,n$ be nonnegative integers such that $m+n\geq 1$. From Theorem \ref{relation}, we get that $a^k$ is left dual core invertible if and only if $a^m$ is left dual $(a^k,a^n)$-core invertible, for some positive integer $k$.

According to Theorem \ref{relation}, we can establish the relation between left dual pseudo core  inverses and left dual $v$-core inverses.
 \begin{prop}
     Let $R$ be a ring, $a\in R$ and $m$ be a nonnegative integer, then $a$ is left dual pseudo core invertible if and only if $a^k$ is left dual $a^m$-core invertible for some positive integer $k$. In this case, $a_{l\bccore{D}}=(a^k)_{l,a^m,\bccore{\#}} a^{k+m-1}$ and $(a^k)_{l,a^m,\bccore{\#}}=(a_{l,\bccore{D}})^{k+m}$.
 \end{prop}

 \begin{proof}
     It is clear from Theorem \ref{relation}.
 \end{proof}
As a result of the proposition, if we take $k=1$, then we have following.
\begin{cor}
   Let $R$ be a ring, $a\in R$ and $m$ be a nonnegative integer. Then $a$ is left dual core invertible if and only if a is left dual $a^m$-core invertible.
\end{cor}

\begin{cor}
    Let $R$ be a ring and $a \in R$. Then following statements are equivalent:
    \begin{enumerate}
        \item $a\in R_{l,\bccore{\#}}$.
        \item $a$ is left dual $(a,a)$-core invertible.
        \item $a$ is left dual $(a,1)$-core invertible.
        \item $1$ is left dual $(a,a)$-core invertible.
        \item $a$ is left dual $a$-core invertible.
        \item $a$ is left dual $1$-core invertible.
        \item $a$ is left $(a^*,a)$- invertible.
    \end{enumerate}
\end{cor}

\begin{cor}
    Let $A\in M_n(\mathbb{C})$ with $I(A)=k$. Then 
    \begin{enumerate}
        \item $A^m$ is left dual $(A^k,A^n)$-core invertible for any nonnegative integers $m,n$ satisfying $m+n\geq 1$. In this case, $(A_{l,\bccore{D}})^{m+k}$ is a left dual $(A^k,A^n)$-core inverse of $A^m$.
        \item $A^k$ is left dual $A^m$-core invertible, for any nonnegative integer $m$. In this case, $(A_{l,\bccore{D}})^{k+m}$ is a left dual $(A^k,A^n)$-core inverse of $A^m$.
    \end{enumerate}
\end{cor}

\begin{cor}
    Let $R$ be a ring and  $a,v\in R$. The following statements are equivalent:
    \begin{enumerate}
        \item $a\in R_{l,v,\bccore{\#}}$.
        \item There exists some $x\in Ra$ such that $axva=a$, $(xav)^*=xav$ and $xvax=x$.
        \item There exists some $x\in R$ such that $axva=a$, $Rx\subseteq Ra$ and $a^*R\subseteq xR$.
        \item There exists some $x\in R$ such that $axva=a$, $Rx\subseteq Ra$ and $axva=a$, $l(x)=l(a^*)$ and $Ra\subseteq Ra$.
        \item There exists some $x\in R$ such that $axva=a$, $Rx\subseteq Ra$ and $xR\subseteq a^*R$.
        \item There exists some $x\in R$ such that $axva=a$, $Rx\subseteq Ra$ and $l(a^*)\subseteq l(x)$.
        \item There exists a projection $q\in R$ and an idempotent $p\in R$ such that $Ra\subseteq Rq\subseteq Rva$, $Rp\subseteq Ra$ and $vaR\subseteq pR$.
        In this case $a_{l,v,\bccore{\#}}=q(va)^-p$ for any $(va)^-\in (va)\{1\}$.
    \end{enumerate}
\end{cor}

As a consequence of above corollary, we have the following result.

\begin{cor}
Let $R$ be a ring and $a\in R$. Then the following statements are equivalent:
    \begin{enumerate}
        \item $a\in R_{l,\bccore{\#}}$.
        \item There exists some $x\in Ra$ such that $axa=a$, $(xa)^*=xa$ and $xax=x$.
        \item There exists some $x\in R$ such that $axa=a$, $Rx\subseteq Ra$ and $a^*R\subseteq xR$.
        \item There exists some $x\in R$ such that $axa=a$, $Rx\subseteq Ra$ and $axa=a$, $l(x)=l(a^*)$ and $Ra\subseteq Ra$.
        \item There exists some $x\in R$ such that $axa=a$, $Rx\subseteq Ra$ and $xR\subseteq a^*R$.
        \item There exists some $x\in R$ such that $axa=a$, $Rx\subseteq Ra$ and $l(a^*)\subseteq l(x)$.
        \item There exists a projection $q\in R$ and an idempotent $p\in R$ such that $Ra\subseteq Rq\subseteq Ra$, $Rp\subseteq Ra$ and $aR\subseteq pR$.
        In this case $a_{l,\bccore{\#}}=q(a)^-p$ for any $(a)^-\in (a)\{1\}$.
    \end{enumerate}
\end{cor}

\begin{thm}
    Let $R$ be a ring and $a\in R$. Then the following statements are equivalent.

    \begin{enumerate}
        \item $a$ is Moore-Penrose invertible.
        \item $a$ is left dual $a^*$-core invertible.
        \item $a^*$ is left dual $a$-core invertible.
        \item $a$ is left dual $(a^*,a^*)$-core invertible.
        \item $a$ is left $(a^*,a^*)$ invertible.
        \item $a^*$ is left dual $(a,a)$-core invertible.
        \item $a^*$ is left $(a,a)$ invertible.
    \end{enumerate}
\end{thm}

\begin{proof}
    (1) $\Leftrightarrow$ (2) $\Leftrightarrow$ (3) by \cite[ Theorem 2.21]{Mosic}.

    (1) $\Rightarrow$ (4) Let $a\in R^{\dag}$. Then we have $a^* \in R^{\{1,4\}}$. Again, since $a\in R^{\dag}$, $a$ is left invertible along $a^*$ by \cite[Corollary 2.21]{Patricio}. Also, $a^*\in Ra^*aa^*$ by \cite[Theorem 2.3]{Patricio}. Hence, $a\in R_l^{(a^*,a^*)}$. Consequently, $a\in _{l,\bccore{\#}\!(a^*,a^*)}$ by Theorem \ref{hehehe}.

    (4) $\Rightarrow$ (5) It is clear from Theorem \ref{hehehe}.

    (5) $\Rightarrow$ (1) Let $a\in R_l^{(a^*,a^*)}$. Then $a^*\in Ra^*aa^*$. Hence, we get $a\in aa^*aR$ and by \cite[Lemma 2.20]{Mosic}, $a\in R^{\dag}$.

    (1) $\Leftrightarrow$ (6) $\Leftrightarrow$ (7) If we consider the proof of (1) $\Leftrightarrow$ (4) $\Leftrightarrow$ (5) and the fact $a\in R^{\dag}$ if and only  $a^*\in R^{\dag}$, then the proof is clear.
\end{proof}

\begin{thm}\label{serhat}
    Let $R$ be a ring and $a,b,c \in R $. Then $a$ is left dual $(b,c)$-core invertible if and only if $ab$ is left $(b^*,c)$ invertible. In this case, the left dual $(b,c)$-core inverse of $a$ coincides with the left $(b^*,c)$ inverse of $ab$.
\end{thm}

\begin{proof}
    Let $a\in R_{l,\bccore{\#}\!(b,c)} $ and $a_{l,\bccore{\#}\!(b,c)}=x$. Then $x\in Rc$, $bxab=b$ and $(xab)^*=xab$. Thus, $b^*=(bxab)^*=(xab)^*b^*=xabb^*\in Rcabb^*$ as required.

    Conversely,  let $x=(ab)_l^{(b^*,c)} $. Then $x\in Rc$ and $xabb^*=b^*$. So since $b=b(xab)^*$, $xab=xab(xab)^*$. Hence $xab=(xab)^*$. Finally $bxab=b(xab)^*=b$ as desired.
\end{proof}

Recall that an element $a\in R$ is \textit{strongly left $(b,c)$ invertible} if $b\in Rcab$ and $cab$ is regular, or equivalently if there exists $x\in R$ such that $xax=x$, $xR=bR$, and $Rx\subseteq Rc$. In which case any such $x$ will be called strongly left $(b,c)$ inverse of $a$.  Clearly, if $x$ is strongly left $(b,c)$ inverse of $a$, then so is $xax$. It is clear that every strongly left $(b,c)$ invertible element is strongly left $(b,c)$ invertible. Moreover, every strongly left $(b,c)$ inverse of $a$ is a left $(b,c)$ inverse of $a$.

Now we have the following theorem.

\begin{thm}
    Let $R$ be a ring and $a,b,c\in R$. Then $a$ is left dual $(b,c)$-core invertible if and only if $ab$ is strongly left $(b^*,c)$ invertible.
\end{thm}

\begin{proof}
Assume that  $ab$ is strongly left $(b^*,c)$ invertible. So $ab$ is left $(b^*,c)$ invertible. By Theorem \ref{serhat}, $a$ is left dual $(b,c)$-core invertible. For the other direction, let $a\in R_{l,\bccore{\#}\!(b,c)}$ and $a_{l,\bccore{\#}\!(b,c)}=y$. Then $y\in Rc$, $byab=b$ and $(yab)^*=yab$, whence $y=tc$ for some $t\in R$. Consequently, $b^*=yabb^*=tcabb^*\in Rc(ab)b^*$. Now, $c(ab)b^*=c(ab)(tcabb^*)=ca(tcabb^*)^*(tcabb^*)=cabb^*(tca)^*tcabb^*$. So, $cabb^*$ is regular and this completes the proof.
\end{proof}

\begin{cor}
    Let $R$ be a ring and $a\in R$. The following statements are equivalent.

    \begin{enumerate}
        \item $a\in  R_{l,\bccore{\#}}$.
        \item $a^2$ is strongly left $(a^*,a)$ invertible.
        \item $a^2$ is strongly left $(a^*,1)$ invertible.
        \item $a$ is strongly left $(a^*,a)$ invertible.
    \end{enumerate}
\end{cor}

\begin{cor}
    Let $R$ be a ring and $a\in R$. The following statements are equivalent.

    \begin{enumerate}
        \item $a$ is Moore-Penrose invertible.
        \item $a^*a$ is strongly left $(a^*,a)$ invertible.
        \item $aa^*$ is strongly left $(a,a^*)$ invertible.
    \end{enumerate}
\end{cor}

The next corollary is a direct consequence of \cite[Theorem 2.21]{Mosic} and Theorems \ref{relation} and \ref{dec}.

\begin{cor}
    Let $R$ be a ring and $a\in R^{\dag}$. Then $aa^*=a_1+a_2$ where 

    \begin{enumerate}
        \item $a_1\in R_{l,\bccore{\#}} $,
        \item $a_2^2=0$,
        \item $a_2^*a_1=0=a_1a_2$. In addition, $a^*a\in R_{l,\bccore{\#}}$ with a left dual core inverse $a^{\dag}{a^{\dag}}^*$.
    \end{enumerate}
\end{cor}

\begin{cor}
    Let $R$ be a ring and $a\in R_{l,\bccore{D}}$ with $I(a)=k$. Then $a^k=a_1+a_2$ where

    \begin{enumerate}
        \item $a_1\in R_{l,\bccore{\#}}$,
        \item $a_2^2=0$,
        \item $a_2^*a_1=0=a_1a_2$. In addition, $a_{l,\bccore{D}}a^{k+1}\in R_{l,\bccore{\#}}$ with a left dual core inverse $(a_{l,\bccore{D}})^k$.
    \end{enumerate}
\end{cor}

\begin{thm}\label{son}
Let $R$ be a ring and $a,d\in R_l^{(b,c)}$. Then the following hold.

\begin{enumerate}
    \item $a_l^{(b,c)}=d_l^{(b,c)}da_l^{(b,c)}$,
    \item $d_l^{(b,c)}=a_l^{(b,c)}ad_l^{(b,c)}$.
\end{enumerate}
\end{thm}

\begin{proof}
    Let $a,d\in R_l^{(b,c)}$, $x=a_l^{(b,c)}$ and $y=d_l^{(b,c)}$. Then we have $x,y\in Rc$ and $b=xab=ydb$.
    
   (1) Let $z=ydx$. Now, $z\in Rc$ since $x\in Rc$. Also, $zab=ydxab=ydb=b$. This implies that $ydx$ is one of the $(b,c)$ inverses of $a$.

   (2) Let $t=xay$. Since $y\in Rc$, we have $t\in Rc$. Also, $tdb=xaydb=xab=b$. This implies that $xay$ is one of the $(b,c)$ inverses of $a$, as desired.
\end{proof}

\begin{thm}
  Let $R$ be a ring and $a,d \in R_{l,\bccore{\#}\!(b,c)}$. Then we have the following.
  \begin{enumerate}
      \item $ a_{l,\bccore{\#}\!(b,c)}=d_{l,\bccore{\#}\!(b,c)}da_l^{(b,c)}$.
      \item $d_{l,\bccore{\#}\!(b,c)}=a_{l,\bccore{\#}\!(b,c)}ad_l^{(b,c)}$.
  \end{enumerate}
\end{thm}

\begin{proof}
    Let $a,d\in R_{l,\bccore{\#}\!(b,c)}$. Then by Theorem \ref{hehehe}, $a_{l,\bccore{\#}\!(b,c)}=b^{(1,4)}a_l^{(b,c)}$ and $d_{l,\bccore{\#}\!(b,c)}=b^{(1,4)}d_l^{(b,c)}$. Also we have $a_l^{(b,c)}=d_l^{(b,c)}da_l^{(b,c)}$ and $d_l^{(b,c)}=a_l^{(b,c)}ad_l^{(b,c)}$ by Theorem \ref{son}. Hence,

    (1) $a_{l,\bccore{\#}\!(b,c)} =b^{(1,4)}a_l^{(b,c)}=b^{(1,4)}d_l^{(b,c)}da_l^{(b,c)}=d_{l,\bccore{\#}\!(b,c)}da_l^{(b,c)}$.

    (2) $d_{l,\bccore{\#}\!(b,c)}= b^{(1,4)}d_l^{(b,c)}=b^{(1,4)}a_l^{(b,c)}ad_l^{(b,c)}=a_{l,\bccore{\#}\!(b,c)}ad_l^{(b,c)}$.
\end{proof}

\begin{thm}
    Let $R$ be a ring and $a,b,c\in R$. Then following are equivalent.
    \begin{enumerate}
        \item $a\in R_{l,\bccore{\#}\!(b,c)}\cap R_{r(b,c)}^{\bccore{\#}\!}$.
        \item $a\in R_{\bccore{\#}\!(b,c)}\cap R_{(b,c)}^{\bccore{\#}\!}$.
        \item $cab\in R^{\dag}$.
        \item $b\in (cab)^*R$ and $c \in R(cab)^*c$.
        \item $R=R(cab)^*\oplus l(c)=(cab)^*R\oplus r(b)$.
        \item $R=R(cab)^*+ l(c)=(cab)^*R+ r(b)$.
    \end{enumerate}
\end{thm}

\begin{proof}
    It is sufficient to show (1) $\Rightarrow$ (2). Because (2) $\Rightarrow$ (1) is obvious and (2) $\Leftrightarrow$ (3) $\Leftrightarrow$ (4) $\Leftrightarrow$ (5) $\Leftrightarrow$ (6) follows from \cite{Zhu}. Now let $a\in R_{l,\bccore{\#}\!(b,c)}\cap R_{r,(b,c)}^{\bccore{\#}\!}$. Then by Theorem \ref{hehehe} and \cite{Dong}, we have $a\in R_l^{(b,c)}\cap R_r^{(b,c)}$, $b\in R^{(1,4)}$ and $c\in R^{(1,3)}$. Since $a\in R_l^{(b,c)}\cap R_r^{(b,c)}$, $a\in R^{(b,c)}$. Hence, $a\in R_{\bccore{\#}\!(b,c)}\cap R_{(b,c)}^{\bccore{\#}\!}$ as desired.
\end{proof}


\begin{thebibliography}{99}

\bibitem{D2}M.P. Drazin, {\it A class of outer generalized inverses}, Linear Algebra Appl., 436(7) (2012), 1909-1923.







\bibitem{Zhang} H. H. Zhu, X.X. Zhang and J.L. Chen, {\it Generalized inverses of a factorization in a ring with involution}, Linear Algebra Appl., 472(2015), 142-150.

\bibitem{Zhu} H.H. Zhu, {\it The (b c)-core inverse and its dual in semigroups with involution} , J. Pure Appl. Algebra, 228(4) (2024), 107526.

\bibitem{Drazin} M.P. Drazin, {\it Left and right generalized inverses} , Linear Algebra Appl., 510 (2016), 64-78.

\bibitem{Mosic} H.H. Zhu, L.Y. Wu and D. Mosic, {\it One-sided w-core inverses in rings with
 involution}, Linear Multilinear Algebra, 71(4) (2023), 528-544.

\bibitem{Patricio} Zhu HH, Chen JL, Patrício P. {\it Further results on the inverse along an element in semigroups and rings}, Linear Multilinear Algebra, 64(3) (2016), 393–403.

\bibitem{Dong} Dong B, Jin Taohua and Zhu HH. {\it Right (b,c)-core inverses in rings}, Submitted for publication.

\bibitem{Mary} X. Mary, {\it On generalized inverses and Green's relations}, Linear Algebra Appl., 434(8) (2011), 1836-1844.

\bibitem{D1}  M.P. Drazin, {\it Pseudo-inverses in associative rings and semigroups}, Amer. Math. Monthly, 65 (1958), 506-514.

\bibitem{Bentez} J. Benítez and E. Boasso, {\it The inverse along an element in rings}, Electron. J. Linear Algebra, 31 (2016), 572-592.

\bibitem{Boasso}  J. Benítez and E. Boasso, {\it The inverse along an element in rings with an involution, Banach algebras and C-algebras}, Linear Multilinear Algebra, 65(2) (2017), 284-299.

\bibitem{Ilic}  D.S. Cvetkovic-Ilic, Y.M. Wei, Algebraic Properties of Generalized Inverses, Series: Dev. Math., Vol. 52, Springer, 2017.

\bibitem{Marypat} X. Mary and P. Patricio, {\it Generalized inverses modulo H in semigroups and rings} , Linear Multilinear Algebra, 61(8) (2013), 1130-1135.

\bibitem{Penrose} R. Penrose, {\it A generalized inverse for matrices} , Proc. Camb. Phil. Soc.,
 51 (1955), 406-413.

\bibitem{Baksalary} O.M. Baksalary and G. Trenkler, {\it  Core inverse of matrices}, Linear Multilinear Algebra,  58(5-6) (2010), 681-697.

\bibitem{Mohana} K. Manjunatha Prasad and K.S. Mohana, 
{\it Core-EP inverse} , Linear Multilinear Algebra 62(6)(2014), 792-802.

\end{thebibliography}
\end{document}